\documentclass[a4paper,12pt]{article}
\usepackage{amsmath,amscd,amssymb,epsfig}
\usepackage[francais]{babel}

\def\a{\alpha}

\def\wb{\overline}
\def\wt{\widetilde}
\def\Mod{\hbox{Mod}}

\def\G{{\cal{G}}}

\def\T{{\cal{T}}}

\def\spin{\hbox{\tiny{spin}}}
\def\ZZ{\mathbb{Z}}  
\def\QQ{\mathbb{Q}}

\def\CC{\mathbb{C}}
\def\NN{\mathbb{N}}
\def\go{\longrightarrow}
\def\im{\mapsto}

\def\A{{\cal{A}}}
\def\B{{\cal{B}}}

\def\C{{\cal{C}}}
\def\D{{\cal{D}}}

\def\L{\Lambda}

\def\S{\mathfrak{S}}

\def\sll{\mathfrak{sl}}
\def\psl{\mathfrak{psl}}

\def\osp{\mathfrak{osp}}
\def\sp{\mathfrak{sp}}
\def\so{\mathfrak{so}}
\def\spin{\hbox{{\footnotesize spin}}}

\def\gg{{\mathfrak{g}}}
\def\ff{{\mathfrak{f}}}
\def\ee{{\mathfrak{e}}}

\newtheorem{theo}{\bf Th{\'e}or{\`e}me}
\newtheorem{prop}[theo]{\bf Proposition}
\newtheorem{propri}[theo]{\bf Propri{\'e}t{\'e}}

\newenvironment{dem}[1][\bf D{\'e}monstration]
  {\begin{trivlist}\item{\em #1.}\ }
  {\hfill$\square$\par\end{trivlist}}
\newenvironment{dem...}[1][\bf D{\'e}monstration]
  {\begin{trivlist}\item{\em #1.}\ }
  {\par\end{trivlist}}

\begin{document}
\begin{center}
{\large{\bf{Invariant d'entrelacs associ{\'e} {\`a} la repr{\'e}sentation des spineurs de $\so_7$\\}}}
\vspace{0.3cm}
{\large{\bf\it{Link invariant for the spinor representation of $\so_7$\\}}}
\vspace{0.5cm}
{\bf{Bertrand Patureau-Mirand\\}}
{\sl \small
Centre de Math{\'e}matiques de Jussieu,
Universit{\'e} Paris 7 Denis Diderot \\
Case Postale 7012, 
2, place Jussieu,
F-75251 PARIS CEDEX 05 \\
E-mail : patureau@math.jussieu.fr \\
URL : http://www.math.jussieu.fr/{$\sim$}patureau\\}
\end{center}
\begin{abstract}
  La fonction de poids associ{\'e}e {\`a} la repr{\'e}sentation des spineurs de $\so_7$
  compos{\'e}e avec l'invariant universel de Vassiliev-Kontsevich fournit un
  invariant d'entrelacs num{\'e}rique {\`a} valeur dans les s{\'e}ries formelles en une
  variable. Je calcule des relations ``skeins'' pour cet invariant et je donne
  un algorithme pour le calculer. Cet invariant est en fait {\`a} valeur dans
  l'anneau $\ZZ[W,W^{-1}]$ des polyn{\^o}mes de Laurent en une variable {\`a}
  coefficients entiers. 
\end{abstract}
\begin{abstract}
  Pulling back the weight system associated with the spinor representation of
  the Lie algebra $\so_7$ by the universal Vassiliev-Kontsevich invariant yields
  a numerical link invariant with values in formal power series. Computing
  some skein relations satisfied by this invariant, I derive a recursive
  algorithm for its evaluation. The values of this invariant belong to the
  ring $\ZZ[W,W^{-1}]$ of Laurent polynomials in one variable. 
\end{abstract}
\section*{Introduction}
Cet article est tir{\'e} de mon travail en th{\`e}se.\\ 
Pour chaque superalg{\`e}bre $L$ de la famille exceptionnelle {\'e}tendue (famille
comprenant les cinq alg{\`e}bres exceptionnelles $\gg_2$, $\ff_4$, $\ee_6$, $\ee_7$
et $\ee_8$ auxquelles on peut ajouter les deux superalg{\`e}bres de Lie $\gg(3)$ et
$\ff(4)$, les superalg{\`e}bres $\psl(E_0)$, $\sll(E_2)$, $\sll(E_3)$,
$\osp(F_{-1})$, $\osp(F_8)$ o{\`u} $E_n$ (respectivement $F_n$) d{\'e}signe un
superespace vectoriel de superdimension $n$ (respectivement muni d'une forme
bilin{\'e}aire supersym{\'e}trique non d{\'e}g{\'e}n{\'e}r{\'e}e), il existe une sous-alg{\`e}bre $l$ de
$L$ et une repr{\'e}sentation de cette alg{\`e}bre $e$ telle que si $v$ est la
repr{\'e}sentation standard de $\sp_2$ on ait~:
$$L\simeq\sp_2\oplus l\oplus v\otimes e$$ 
comme $(\sp_2\times l)$-module. Les fonctions de poids associ{\'e}es aux couples
$(l,e)$ ont des propri{\'e}t{\'e}s communes qui m'ont permis de calculer certaines
relations ``skeins'' v{\'e}rifi{\'e}es par les invariants d'entrelacs associ{\'e}s. Ce
syst{\`e}me de relations n'est malheureusement pas complet (on ne peut pas en
d{\'e}duire directement un calcul de ces invariants). Cependant, pour le cas
$L=\ff(4)$ o{\`u} l'on a $l=\so_7$ et o{\`u} $e$ est la repr{\'e}sentation des spineurs, la
fonction de poids associ{\'e}e au couple $(\so_7,\spin_7)$ a des propri{\'e}t{\'e}s
suppl{\'e}mentaires desquelles se d{\'e}duisent d'autres relations ``skeins''. On a
alors un syst{\`e}me complet de relations duquel on peut d{\'e}duire un algorithme pour
le calcul de l'invariant d'entrelacs associ{\'e}.\\

L'invariant de Vassiliev des noeuds le plus g{\'e}n{\'e}ral associ{\'e} {\`a} une alg{\`e}bre de Lie
$\gg$ est $Z_\gg$ qui est {\`a} valeur dans les s{\'e}ries formelles {\`a} une variable {\`a}
coefficients dans le centre de l'alg{\`e}bre enveloppante de $\gg$. Pour obtenir un
invariant num{\'e}rique, on {\'e}value la trace de l'action de ces {\'e}l{\'e}ments sur une
repr{\'e}sentation. Il est bien connu que la famille d'invariants obtenue pour
$\gg=\osp(E)$ et pour le choix de la repr{\'e}sentation standard $E$ est donn{\'e}e par
le polyn{\^o}me de Kauffman. Les invariants obtenus pour le choix de repr{\'e}sentations
du groupe $SO(E)$ se d{\'e}duisent des diff{\'e}rents cablages du polyn{\^o}me de
Kauffman. Ce n'est pas le cas des repr{\'e}sentations des spineurs qui sont des
repr{\'e}sentations des alg{\`e}bres $\so_n(\CC)$ mais ne sont pas des repr{\'e}sentations
du groupe $SO_n(\CC)$ (l'existence de ces repr{\'e}sentations est li{\'e}e au fait que
le groupe $SO_n(\CC)$ n'est pas simplement connexe). Bien s{\^u}r, la connaissance
de tous les cablages du polyn{\^o}me de Kauffman d{\'e}termine compl{\`e}tement l'{\'e}l{\'e}ment
$Z_{\so_n}$ et donc l'invariant obtenu pour la repr{\'e}sentation des spineurs de
$\so_n$, mais, en pratique, ce calcul est impossible. L'invariant que je pr{\'e}sente
poss{\`e}de donc un int{\'e}r{\^e}t propre et n'est pas directement une sp{\'e}cialisation d'un
cablage du polyn{\^o}me de Kauffman.

\section{La fonction de poids $\Phi_{\so_7,{\hbox{{\small spin}}}_7}$}
\subsection{D{\'e}finition des diagrammes}
Un $(X_1,X_2)$-diagramme bicolore est un graphe fini $K$,
dont tous les sommets sont trivalents ou univalents, muni des donn{\'e}es
suivantes~:
\begin{enumerate}
\item Pour chaque sommet trivalent $x$ de $K$, un ordre cyclique sur
  l'ensemble des trois ar{\^e}tes orient{\'e}es arrivant en $x$.
\item Un isomorphisme identifiant l'ensemble des sommets univalents {\`a} $X_1\amalg
X_2$.
\item Une application $c$ de l'ensemble des ar{\^e}tes de $K$ vers l'ensemble des
``couleurs'' $\{1,2\}$ telle que $c^{-1}(\{2\})$ forme une courbe non orient{\'e}e de
bord $X_2$.
\end{enumerate}

On peut repr{\'e}senter un $(X_1,X_2)$-diagramme bicolore par un graphe
uni-tri-valent immerg{\'e} dans le plan de mani{\`e}re {\`a} ce que l'ordre cyclique {\`a}
chaque sommet soit donn{\'e} par l'orientation du plan. On repr{\'e}sentera d'un trait
plus {\'e}pais les ar{\^e}tes de la deuxi{\`e}me couleur.\\ 
Chaque sommet trivalent d'un diagramme bicolore est soit dit de ``couleur $1$''
(si les trois ar{\^e}tes qui en sont issues sont de couleur $1$), soit dit ``mixte''
(si en sont issues deux ar{\^e}tes de couleurs $2$ et une ar{\^e}te de couleur $1$). La
couleur d'un sommet univalent est la couleur de l'ar{\^e}te qui en est issue.\\ 
On d{\'e}finit le degr{\'e} d'un $(X_1,X_2)$-diagramme bicolore par
$n-s-\frac12\hbox{cardinal}(X_2)$.

\subsection{D{\'e}finition des modules de diagrammes}
On note $\A(X_1,X_2)$ le $\QQ$-espace vectoriel de base les
$(X_1,X_2)$-diagrammes bicolores quotient{\'e} par les relations $(AS)$, $(IHX)$ et
$(STU)$ ci-dessous~: 
\begin{enumerate}
\item Si deux diagrammes de $\C$ ne diff{\`e}rent que par l'ordre cyclique de l'un
de leurs sommets trivalents, leur somme est nulle (relation dite (AS) pour
antisym{\'e}trie).
  $$\begin{array}{cccc}\put(-10,-10) {\epsfbox{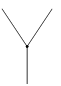}}&+&
    \put(-10,-10) {\epsfbox{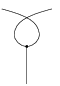}}&\equiv0\end{array}$$
Cette relation est valable quelles que soient les couleurs des trois ar{\^e}tes.
\item la relation (IHX) fait intervenir trois diagrammes de $\C$ qui ne
diff{\`e}rent qu'au voisinage d'une ar{\^e}te~:
  $$\begin{array}{ccccc}\put(-8,-5) {\epsfbox{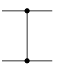}}&\equiv&
    \put(-8,-5) {\epsfbox{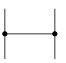}} &
    -&\put(-8,-5) {\epsfbox{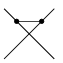}}\end{array}$$
\item La relation (STU) qui est une variation de la relation (IHX) au voisinage
d'un sommet mixte~:
$$\begin{array}{ccccc}\put(-8,-5) {\epsfbox{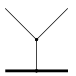}}&\equiv& \put(-8,-5)
{\epsfbox{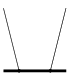}} & -&\put(-8,-5) {\epsfbox{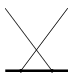}}\end{array}$$
\end{enumerate}
On pose aussi
$$\A(X)=\bigoplus_{X_1\amalg X_2=X} \A(X_1,X_2)$$

\subsection{Cat{\'e}gories de diagrammes}
Soit $\wb\B$ et $\B$ les cat{\'e}gories mono{\"\i}dales $\CC$-lin{\'e}aires d{\'e}finies
par~:
\begin{itemize}
\item Obj$(\wb\B)=$Obj$(\B)=\left\{[n],n\in\NN\right\}$
\item $\wb\B([p],[q])=\A([p]\amalg[q])$\qquad
$\B([p],[q])=\A(\emptyset,[p]\amalg[q])$
\item La composition d'un diagramme de $[p]$ vers $[q]$ avec un diagramme de
$[q]$ vers $[r]$ est nulle si les couleurs des sommets univalents des deux
ensembles $[q]$ ne sont pas identiques. Sinon, elle est donn{\'e}e par la r{\'e}union au
dessus de $[q]$ des deux diagrammes (on les recolle en identifiant les sommets
monovalents de m{\^e}me index des deux ensembles $[q]$).
\item Le produit tensoriel $[p]\otimes [q]$ vaut $[p+q]$ et celui de deux diagrammes
est donn{\'e} par l'image de leur r{\'e}union disjointe par l'isomorphisme de $[p]\amalg
[q]\simeq [p+q]$ obtenu en augmentant de $p$ chaque {\'e}l{\'e}ment de $[q]$.\\ 
\end{itemize}
Les morphismes de la cat{\'e}gorie $\B$ sont donc donn{\'e}s par les combinaisons
lin{\'e}aires de diagrammes dont tous les sommets univalents sont de la deuxi{\`e}me
couleur. On identifiera $\B([p],[q])$ {\`a} un sous-espace de $\wb\B([p],[q])$. Dans
le cas o{\`u} $p=q$, cette inclusion est pour la composition un morphisme d'alg{\`e}bre
non unitaire.

\subsection{Fonctions de poids}
Soit $L$ une (super)alg{\`e}bre de Lie quadratique (munie d'une forme bilin{\'e}aire
invariante non d{\'e}g{\'e}n{\'e}r{\'e}e $<.,.>_L$) et soit $\Omega\in L\otimes L$ l'{\'e}l{\'e}ment de
Casimir associ{\'e}.\\ 
On consid{\`e}re $E$ un $L$-module autodual, muni d'une forme bilin{\'e}aire invariante
non d{\'e}g{\'e}n{\'e}r{\'e}e $<.,.>_E$) et soit $\pi\in L\otimes L$ l'{\'e}l{\'e}ment de Casimir
associ{\'e}.\\

\begin{propri}
Il existe un foncteur mono{\"\i}dal $\CC$-lin{\'e}aire $\wb\Phi_{L,E}$ de $\wb\B$ dans la
cat{\'e}gorie $\Mod_L$ des repr{\'e}sentations de $L$ envoyant $[1]$ sur $L\oplus E$ et
d{\'e}finit de mani{\`e}re unique par ses valeurs prises sur les diagrammes suivant~:
$${\epsfbox{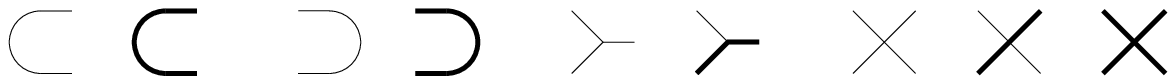}}$$
qui ont pour image respectivement~:
\begin{enumerate}
\item Le Casimir de $L$~: $\Omega\in L^{\otimes 2}\hookrightarrow \Mod_L(\CC,(L\oplus
E)^{\otimes 2})$
\item Le Casimir de $E$~: $\pi\in E^{\otimes 2}\hookrightarrow \Mod_L(\CC,(L\oplus
E)^{\otimes 2})$
\item Le produit scalaire sur $L$ associ{\'e} $\Omega$ de $L^{\otimes 2}\go\CC$ vu
comme morphisme de $\Mod_L((L\oplus E)^{\otimes 2},\CC)$
\item Le produit scalaire sur $E$ associ{\'e} $\pi$ de $E^{\otimes 2}\go\CC$ vu
comme morphisme de $\Mod_L((L\oplus E)^{\otimes 2},\CC)$
\item Le crochet de Lie de $L^{\otimes 2}\go L$ vu comme morphisme de
$\Mod_L((L\oplus E)^{\otimes 2},(L\oplus E))$
\item L'op{\'e}ration externe du $L$-module $E$ de $L\otimes E\go E$
\item Les trois op{\'e}rateurs de sym{\'e}trie ~: $\begin{array}[t]{ccl} X\otimes
Y&\go&Y\otimes X\\ x\otimes y&\im&y\otimes x \end{array}$
\end{enumerate}

De m{\^e}me, on note $\Phi_{L,E}:\B\go \Mod_L$ le foncteur mono{\"\i}dal $\CC$-lin{\'e}aire
envoyant $[1]$ sur $E$ et d{\'e}fini sur les morphismes par~:\\
$\Phi_{L,E}(K)=i^*(p_*(\circ\wb\Phi_{L,E}(K)))$ o{\`u} $p$ d{\'e}signe la projection
canonique de $L\oplus E$ sur $E$ et $i$ d{\'e}signe l'inclusion canonique de $E$
dans $L\oplus E$.
\end{propri}

\subsection{Repr{\'e}sentations de $\so_7$}
On note $\Gamma(p,q,r)$ la repr{\'e}sentation de $L=\so_7$ de plus haut poids
$p\omega_1+q\omega_2+r\omega_3$ o{\`u} $\omega_1$ est le plus haut poids de la
repr{\'e}sentation standard $v$, $\omega_2$ celui de la repr{\'e}sentation adjointe $l$ et
$\omega_3$ celui de la repr{\'e}sentation des spineurs $s$. On a dans la cat{\'e}gorie
des repr{\'e}sentations de $\so_7$~:
$$\Lambda^2s=v\oplus l\qquad S^2s=\CC\oplus \Gamma(0,0,2)$$
$$s\otimes v=s\oplus \Gamma(1,0,1)\qquad s\otimes l=s\oplus \Gamma(1,0,1) \oplus
\Gamma(0,1,1)$$
$$s\otimes \Gamma(0,0,2)=s\oplus \Gamma(1,0,1) \oplus \Gamma(0,1,1)\oplus
\Gamma(0,0,3)$$
$$\Lambda^3s=s\oplus \Gamma(1,0,1)\qquad S^3s=s\oplus \Gamma(0,0,3)$$

Nous prenons comme Casimir de r{\'e}f{\'e}rence quarante fois celui associ{\'e} {\`a} la forme
de Killing. Il agit sur les diff{\'e}rentes repr{\'e}sentations {\'e}voqu{\'e}es par les
scalaires suivants~:\\
$$\begin{array}{|r|c|c|}\hline
\hbox{Repr{\'e}sentation}&\hbox{dimension}&\hbox{action du Casimir}\\ \hline
v=\Gamma(1,0,0)&7&24\\ \hline l=\Gamma(0,1,0)&21&40\\ \hline
s=\Gamma(0,0,1)&8&21\\ \hline \Gamma(0,0,2)&35&48\\ \hline \Gamma(1,0,1)&48&49\\
\hline \Gamma(0,1,1)&112&69\\ \hline \Gamma(0,0,3)&112&81\\ \hline \end{array}$$
On a donc $$\Mod_L(s^{\otimes2},s^{\otimes2})\simeq\CC^4 \hbox{ (comme alg{\`e}bre)}$$
$$\Mod_L(s^{\otimes3},s^{\otimes3})\simeq M_4(\CC)\times M_3(\CC)\times
M_2(\CC)\times\CC$$
Nous appelons $\A'_3$ l'alg{\`e}bre gradu{\'e}e
$\Mod_L(s^{\otimes3},s^{\otimes3})\otimes\CC[\a]\simeq M_4(\CC[\a])\times
M_3(\CC[\a])\times M_2(\CC[\a])\times\CC[\a]$ o{\`u} degr{\'e}$(\a)=1$.\\
La fonction de poids $\Phi_{\so_7,\spin_7}$ associ{\'e}e au choix de $\a$ fois le
Casimir standard induit donc un morphisme d'alg{\`e}bre gradu{\'e}e de $\A_3=\B([3],[3])$
dans l'alg{\`e}bre $\A'_3$.\\
Appelons $i$ l'image dans $\wt\D([1],[1])$ de l'unique $([0],[2])$-diagramme
sans boucle. On note respectivement $r$ et $a$ les deux diagrammes suivants~:\\
$$\epsfbox{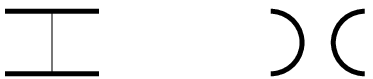}$$
vus comme des morphismes de $\B([2],[2])$ dont la source est {\`a} gauche et le but
{\`a} droite. On appelle aussi $s$ l'op{\'e}rateur de sym{\'e}trie de $\B([2],[2])$

\begin{prop}
La fonction de poids $\Phi_{\so_7,\spin_7}$ annulle les {\'e}l{\'e}ments suivants~:\\
$$\begin{picture}(0,0)%
\includegraphics{Normpex.pstex}%
\end{picture}%
\setlength{\unitlength}{1973sp}%
\begingroup\makeatletter\ifx\SetFigFont\undefined%
\gdef\SetFigFont#1#2#3#4#5{%
  \reset@font\fontsize{#1}{#2pt}%
  \fontfamily{#3}\fontseries{#4}\fontshape{#5}%
  \selectfont}%
\fi\endgroup%
\begin{picture}(11301,678)(1012,-550)
\put(1801,-286){\makebox(0,0)[lb]{\smash{\SetFigFont{12}{14.4}{\rmdefault}{\mddefault}{\updefault}
$-8$
}}}
\put(5551,-286){\makebox(0,0)[b]{\smash{\SetFigFont{12}{14.4}{\rmdefault}{\mddefault}{\updefault}
$-8\a$
}}}
\put(10426,-286){\makebox(0,0)[b]{\smash{\SetFigFont{12}{14.4}{\rmdefault}{\mddefault}{\updefault}
$-40\a$
}}}
\end{picture}
$$ 
$$\begin{picture}(0,0)%
\includegraphics{IHXT.pstex}%
\end{picture}%
\setlength{\unitlength}{1973sp}%
\begingroup\makeatletter\ifx\SetFigFont\undefined%
\gdef\SetFigFont#1#2#3#4#5{%
  \reset@font\fontsize{#1}{#2pt}%
  \fontfamily{#3}\fontseries{#4}\fontshape{#5}%
  \selectfont}%
\fi\endgroup%
\begin{picture}(7138,1524)(557,-898)
\put(1501,-136){\makebox(0,0)[b]{\smash{\SetFigFont{12}{14.4}{\rmdefault}{\mddefault}{\updefault}$+$}}}
\put(4801,-136){\makebox(0,0)[b]{\smash{\SetFigFont{12}{14.4}{\rmdefault}{\mddefault}{\updefault}$+$}}}
\put(6526,-136){\makebox(0,0)[b]{\smash{\SetFigFont{12}{14.4}{\rmdefault}{\mddefault}{\updefault}$+6\a$}}}
\put(3001,-136){\makebox(0,0)[b]{\smash{\SetFigFont{12}{14.4}{\rmdefault}{\mddefault}{\updefault}$-3\a$}}}
\end{picture}
$$
$$\begin{picture}(0,0)%
\includegraphics{relspin7.pstex}%
\end{picture}%
\setlength{\unitlength}{1973sp}%
\begingroup\makeatletter\ifx\SetFigFont\undefined%
\gdef\SetFigFont#1#2#3#4#5{%
  \reset@font\fontsize{#1}{#2pt}%
  \fontfamily{#3}\fontseries{#4}\fontshape{#5}%
  \selectfont}%
\fi\endgroup%
\begin{picture}(13138,3624)(557,-3073)
\put(8176,-61){\makebox(0,0)[b]{\smash{\SetFigFont{12}{14.4}{\rmdefault}{\mddefault}{\updefault}
$\L^2$
}}}
\put(3826,-61){\makebox(0,0)[b]{\smash{\SetFigFont{12}{14.4}{\rmdefault}{\mddefault}{\updefault}
$\L^2$
}}}
\put(11026,-61){\makebox(0,0)[b]{\smash{\SetFigFont{12}{14.4}{\rmdefault}{\mddefault}{\updefault}
$\L^2$
}}}
\put(12976,-61){\makebox(0,0)[b]{\smash{\SetFigFont{12}{14.4}{\rmdefault}{\mddefault}{\updefault}
$\L^2$
}}}
\put(3826,-2386){\makebox(0,0)[b]{\smash{\SetFigFont{12}{14.4}{\rmdefault}{\mddefault}{\updefault}
$\L^3$
}}}
\put(6676,-2386){\makebox(0,0)[b]{\smash{\SetFigFont{12}{14.4}{\rmdefault}{\mddefault}{\updefault}
$\L^3$
}}}
\put(9826,-2386){\makebox(0,0)[b]{\smash{\SetFigFont{12}{14.4}{\rmdefault}{\mddefault}{\updefault}
$\L^3$
}}}
\put(2926,-2311){\makebox(0,0)[rb]{\smash{\SetFigFont{12}{14.4}{\rmdefault}{\mddefault}{\updefault}
$+9\a$
}}}
\put(8926,-2311){\makebox(0,0)[rb]{\smash{\SetFigFont{12}{14.4}{\rmdefault}{\mddefault}{\updefault}
$-60\a^2$
}}}
\put(3076,-286){\makebox(0,0)[rb]{\smash{\SetFigFont{12}{14.4}{\rmdefault}{\mddefault}{\updefault}
$-6\a$
}}}
\put(6526,-211){\makebox(0,0)[rb]{\smash{\SetFigFont{12}{14.4}{\rmdefault}{\mddefault}{\updefault}
$-6\a$
}}}
\put(10276,-211){\makebox(0,0)[rb]{\smash{\SetFigFont{12}{14.4}{\rmdefault}{\mddefault}{\updefault}
$+36\a^2$
}}}
\end{picture}
$$ 
Pour ce dernier {\'e}l{\'e}ment, les ``boites'' $\L^2$ (respectivement $\L^3$) d{\'e}signent
la combinaison lin{\'e}aire de diagrammes de $\B([2],[2])$ (respectivement
$\B([3],[3])$) donn{\'e}e par
$\frac12\sum_{\sigma\in\S_2}\epsilon(\sigma)\sigma\,(=\frac12(1-s))$
(respectivement $\frac16\sum_{\sigma\in\S_3}\epsilon(\sigma)\sigma$).\\

De plus, cette repr{\'e}sentation de $\A_3$ est {\'e}quivalente {\`a} celle donn{\'e}e par~:\\
$$r\otimes 1\mapsto r'_{1,2}\qquad s\otimes 1\mapsto s'_{1,2}\qquad 1\otimes
s\mapsto s'_{2,3}\quad\hbox{et}\quad a\otimes 1\mapsto a'_{1,2}$$
avec
$$r'_{1,2}=\left(\left[{\begin{array}{cccc} 21\a & 10\a
& 8\a & 0 \\ 0 & \a & - 4\a & 0 \\ 0 & 0
& 9\a & 6\a \\ 0 & 0 & 0 & - 3\a \end{array}}\right],\right.$$
$$\left.\left[ {\begin{array}{ccc} \a & -20\a & - 12\a \\ 0
& 9\a & 6\a \\ 0 & 0 & - 3\a\end{array}}\right] ,\left[ {\begin{array}{cc} - 3\a
& - 2\a \\ 0 & \a\end{array}}
\right],\left[ {\begin{array}{c} -3\a\end{array}} \right]\right)$$
$$s'_{1,2}=\left(\left[ {\begin{array}{rrrr} 1 & 1 & 1 & 1 \\ 0 & -1 & 0 & -1 \\
0 & 0 & -1 & -1 \\ 0 & 0 & 0 & 1\end{array}} \right], \left[ {\begin{array}{rrr}
-1 & 0 & 1 \\ 0 & -1 & -1 \\ 0 & 0 & 1\end{array}}\right], \left[
{\begin{array}{rr} 1 & 1 \\ 0 & -1\end{array}} \right],\left[ {\begin{array}{c} 1
\end{array}} \right]\right)$$
$$s'_{2,3}= \left(\left[ {\begin{array}{rrrr} 1 & 0 & 0 & 0 \\ -1 & -1 & 0 & 0
\\ -1 & 0 & -1 & 0 \\ 1 & 1 & 1 & 1\end{array}} \right], \left[
{\begin{array}{rrr} 1 & 0 & 0 \\ -1 & -1 & 0 \\ 1 & 0 & -1\end{array}} \right]
,\left[ {\begin{array}{rr} -1 & 0 \\ 1 & 1\end{array}}\right],\left[
{\begin{array}{c} 1 \end{array}} \right]\right)$$
$$a'_{1,2}=\left(\left[{\begin{array}{cccc} 8 & 4&
4 & 1 \\ 0 & 0 & 0 & 0 \\ 0 & 0 & 0 & 0 \\ 0 & 0 & 0 & 0
\end{array}}\right],\left [ {\begin{array}{ccc} 0 & 0 & 0 \\ 0 & 0 & 0 \\ 0 & 0
& 0 \end{array}}\right] ,\left[ {\begin{array}{cc} 0 & 0 \\ 0 & 0 \end{array}}
\right] ,\left[ {\begin{array}{c} 0 \end{array}}
\right]\right) $$
\end{prop}
\begin{dem}
On se contente de justifier l'existence de ces relations v{\'e}rifi{\'e}es par
$\Phi_{\so_7,\spin_7}$ et la m{\'e}thode utilis{\'e}e pour d{\'e}terminer la
repr{\'e}sentation~:\\
Les trois premi{\`e}res relations peuvent {\^e}tre lues de la mani{\`e}re suivantes~: La
repr{\'e}sentation $\spin_7$ est de dimension $8$, la forme bilin{\`e}aire sur $\so_7$
induite par la trace sur $\spin_7$ vaut $8$ fois notre forme bilin{\'e}aire de
r{\'e}f{\'e}rence et la forme de Killing vaut $40$ fois notre forme bilin{\'e}aire de
r{\'e}f{\'e}rence.\\
L'existence des deux derni{\`e}res relations provient du fait que les espaces
$\Mod_L(S^2L,S^2L)$ et $\Mod_L(\L^3L,\L^3L)$ sont de dimension $2$.\\

On utilise ces relations pour calculer le produit dans le quotient de $\A_3$ par
le noyau de $\Phi_{\so_7,\spin_7}$ ce qui permet {\`a} {\'e}quivalence pr{\`e}s de
d{\'e}terminer la repr{\'e}sentation de $\A_3$ dans $\A'_3$.
\end{dem}

\section{L'invariant universel de Vassiliev-Kontsevich}
Soient $R=\CC[[\a]]$ l'anneau des s{\'e}ries en $\a$ {\`a} coefficients dans $\CC$,
$K=\CC[\a^{-1}][[\a]]$ le corps des s{\'e}ries de Laurent en $\a$ {\`a} coefficients
dans $\CC$.\\
Nous reprenons la notion de ``$q$-tangles'' introduite dans [LM1] en n{\'e}gligeant
l'orientation~:\\
Soit $\T_3$ le mono{\"\i}de engendr{\'e} par les $5$ ``$q$-tangles'' {\'e}l{\'e}mentaires, non
orient{\'e}s, munis des parenth{\`e}sages ((**)*), suivants~:
$$R_{1,2}=\put(5,-15){\epsfbox{B17}}\hspace{1cm},
\quad R_{1,2}^{-1}=\put(5,-15){\epsfbox{B20}}\hspace{1cm}, 
\quad R_{2,3}=\put(5,-15){\epsfbox{B18}}\hspace{1cm}, 
\quad R_{2,3}^{-1}=\put(5,-15){\epsfbox{B19}}\hspace{1cm}, 
\quad A_{1,2}=\put(5,-15){\epsfbox{B33}}\hspace{1cm}$$
Notons que $R_{1,2}$, $R_{2,3}$ et leurs inverses engendrent le groupe d'Artin
$B_3\subset\T_3$ (groupe des tresses {\`a} trois brins).\\
Dans ce qui suit, si $A$ est une $\CC$-alg{\`e}bre, $<A>$ d{\'e}signe l'alg{\`e}bre
compl{\'e}t{\'e}e de $A$ pour la graduation. Pour le choix d'un associateur concentr{\'e} en
degr{\'e} pair ({\'e}l{\'e}ment de $< \A_3>$), on note $Z:\T_3\go<\A_3>$ l'application
induite par l'invariant universel de Vassiliev-Kontsevich (cf [LM1]) et
$Z_{\spin_7}:\T_3\go <\A'_3>$ sa compos{\'e}e avec $\Phi_{\so_7,\spin_7}$. soient
$R'_{i,j}=Z_{\spin_7}(R_{i,j})$ et $A'_{i,j}=Z_{\spin_7}(A_{i,j})$.\\
On note encore $Z_{\spin_7}$ l'invariant d'entrelacs en bande non orient{\'e}s {\`a}
valeurs dans $R$.\\

\section{Relations ``skein'' pour $Z_{\spin_7}$}
Remarquons que par restriction, $Z_{\spin_7}$ induit des repr{\'e}sentations
lin{\'e}aires $ \zeta_i:B_3\go\hbox{GL}_{i}(R)$ pour $i\in\{1,2,3,4\}$.

\begin{prop}
Les repr{\'e}sentations $\zeta_i$ sont simples.
\end{prop}
\begin{dem}
La d{\'e}monstration repose sur la connaissance de
$R'_{1,2}=\exp(-\frac\a2r'_{1,2})s'_{1,2}$ et de $R'_{2,3}$ en degr{\'e} inf{\'e}rieur
ou {\'e}gal {\`a} $1$~: en effet, notre associateur est {\'e}gal {\`a} $1$ en degr{\'e} inf{\'e}rieur ou
{\'e}gal {\`a} $1$ et donc $R'_{2,3}\equiv(1-\frac\a2r'_{2,3})s'_{2,3}$ modulo $(\a^2)$.
On peut d'abord remarquer que ${R'_{1,2}}^2\equiv{R'_{2,3}}^2\equiv1$ modulo
$\a$. Ainsi ces repr{\'e}sentations de $B_3$ induisent des repr{\'e}sentations de
$\S_3$. On v{\'e}rifie ais{\'e}ment que $\zeta_2$ est la repr{\'e}sentation simple standard
(associ{\'e}e au tableau de Young $[2,1]$).\\
Pour $\zeta_3$ et $\zeta_4$, il existe une base $(v_1,v_2,v_3,v_4)$ de ${R}^4$
(respectivement $(v_1,v_2,v_3)$ de ${R}^3$) dans laquelle $\zeta_i(R_{1,2})$ est
diagonale et a des valeurs propres distinctes. Ainsi, si un sous espace propre
$V$ de $R^4$ (resp. $R^3$) est stable par $\zeta_i(B_3)$, il se d{\'e}compose en la
somme directe de droites $\bigoplus_{j\in J}R.v_j$. Maintenant, on calcule
$\zeta_i(R_{2,3})$ modulo $\a$ dans les bases $(v_i)$~:
$$\frac{1}{8}\left[ {\begin{array}{rrrr}1&9&3&-15\\3&-5&1&-5\\7&7&1&35\\-1&-1&1&3\end{array}}\right]\qquad\frac {1}{4}\left[{\begin{array}{rrr} -2&-2&10\\-1&-3&-5\\1&-1&1\end{array}} \right] $$
Et l'existence d'un sous-$B_3$-module de type $\wb V$ appara{\^\i}t clairement
impossible.
\end{dem}

Introduisons l'{\'e}l{\'e}ment $W=e^{-\frac12\a}\in R$
\begin{prop} \label{deco matrice}
La repr{\'e}sentation $Z_{\spin_7}$ de $B_3$ dans $\A'_3$ est {\'e}quivalente {\`a} celle
donn{\'e}e par
$$R'_{1,2}=\left(\left[\begin{array}{cccc}W^{21}&-W(1-W^4+W^8)&-W^9(1-W^4+W^8)&W^{-3}\\
0&-W&-(1-W^4)W^9&W^{-3}\\0&0&-W^9&W^{-3}\\0&0&0&W^{-3}\end{array}\right],\right.$$
$$\left.\left[{\begin{array}{ccc}-W&-W^{-11}(W^{20}-1)&-W^9\\ 0&-W^9&-W^9\\
0&0&W^{-3}\end{array}}\right],\left[{\begin{array}{cc}-W&-W\\
0&W^{-3}\end{array}}\right],\left[W^{-3}\right]\right)$$

$$R'_{2,3}=\left(\left[\begin{array}{cccc}W^{-3}&0&0&0\\ W^9&-W^9&0&0\\
W^{-3}&W^{-3}(W^4-1)&-W&0\\
W^9&-W^9(W^8-W^4+1)&-W^{13}(W^8-W^4+1)&W^{21}\end{array}\right],\right.$$
$$\left.\left[{\begin{array}{ccc}W^{-3}&0&0\\ W^9&-W^9&0\\
-W^9&W^{-11}(W^{20}-1)&-W\end{array}}\right],\left[{\begin{array}{cc}W^{-3}&0\\
-W^{-3}&-W\end{array}}\right],\left[W^{-3}\right]\right)$$

\noindent De plus, $Z_{\spin_7}(A_{1,2})$ est donn{\'e} dans la m{\^e}me base par~: 
$$\begin{array}{ll}A'_{1,2}&=\left(\frac\Delta{(W^{20}+1)(W^4+1)}\right.\\ \\
&\left[\begin{array}{cccc}
(W^{20}+1)(W^4+1)&-(W^{12}+1)&-W^4(W^{16}+1)&\frac{W^{12}}{W^{12}+1}\\ 0&0&0&0\\
0&0&0&0\\ 0&0&0&0\end{array}\right],\\
&\left.\left[{\begin{array}{ccc}0&0&0\\0&0&0\\0&0&0 \end{array}}\right],
\left[{\begin{array}{cc}0&0\\0&0\end{array}}\right],\left[0\right]\right)\end{array}$$
\end{prop}
\begin{dem}
Justifions d'abord la forme de $R'_{1,2}$ et $R'_{2,3}$~:\\
Dans [TW], il est {\'e}tabli qu'il existe {\`a} conjugaison pr{\`e}s, une unique
repr{\'e}senta\-tion simple du groupe $B_3$ dans $GL_n(K)$, o{\`u} $K$ est un corps
alg{\'e}briquement clos, pour $n\leq3$, une fois fix{\'e}e la liste
$(\lambda_1,\ldots,\lambda_n)\in K^n$ des valeurs propres des tresses
{\'e}l{\'e}mentaires. Ceci donne le r{\'e}sultat annonc{\'e} pour les $3$ derniers facteurs.\\
Pour $n=4$, la repr{\'e}sentation de $B_3$ est d{\'e}termin{\'e}e par le choix d'une racine
carr{\'e}e $\epsilon W^{-4}$ de
$\frac{\lambda_2\lambda_3}{\lambda_1\lambda_4}=W^{-8}$ (cf [TW]). Or nous savons
que la repr{\'e}sentation $\zeta_4$ de $B_3$ est en fait {\`a} valeur dans $M_4(R)$ et
qu'elle envoie $R_{1,2}$ sur une matrice diagonalisable dans $M_4(R)$. Les
espaces propres de l'image de $R_{1,2}$ par les repr{\'e}sentations correspondantes
{\`a} $\epsilon=1$ et $\epsilon=-1$ sont de la forme $(R v_i)_{i=1\ldots4}$ et le
calcul donne pour $\epsilon=1$, $\det(v_1,v_2,v_3,v_4)\in\alpha^2R$, donc
$R_{1,2}$ n'est alors pas diagonalisable dans $M_4(R)$. La repr{\'e}sentation
$\zeta_4$ est donc {\'e}quivalente {\`a} celle obtenue pour $\epsilon=-1$.

D'autre part, $Z_{\spin_7}(A_{1,2})$ est envoy{\'e} sur $\varphi a'_{1,2}=8\varphi
p_E$ o{\`u} $p_E$ est le projecteur sur $E^4\subset s^{\otimes 3}$ qui est l'espace
propre de $R'_{1,2}$ associ{\'e} {\`a} la valeur propre $W^{21}$ et $\varphi\in
R^*$. Or, on a dans $\T_3$ l'identit{\'e} $A_{1,2}R_{1,2}R_{2,3}A_{1,2}=A_{1,2}$ qui
donne\\ $a'_{1,2}R'_{1,2}R'_{2,3}a'_{1,2}=\varphi^{-1}a'_{1,2}$ ce qui
permet de d{\'e}terminer $\Delta=8\varphi$~:
$$\Delta=W^{-18}(W^{20}+1)(W^{12}+1)(W^4+1)$$ 
\end{dem}

\begin{figure}[p]
$$\begin{array}{|l|}
\hline\\ 
 \put(5,-15){\epsfbox{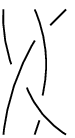}}\hspace{1cm} 
 -\put(5,-15){\epsfbox{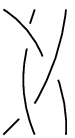}}\hspace{1cm}= \\ \\ 
{W}^{-1}({W}^{4}-1)\\
\hfill\left(\put(5,-15){\epsfbox{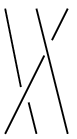}}\hspace{1cm}
 -\put(5,-15){\epsfbox{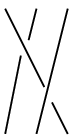}}\hspace{1cm}
 +\put(5,-15){\epsfbox{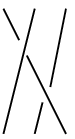}}\hspace{1cm}
 -\put(5,-15){\epsfbox{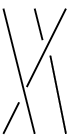}}\hspace{1cm}
 +\put(5,-15){\epsfbox{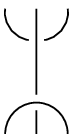}}\hspace{1cm}
 -\put(5,-15){\epsfbox{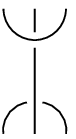}}\hspace{1cm} \right) \\ \\ 
+{W}^{-10}({W}^{4}-1)({W}^{12}-{W}^{8}-1)\\
\hfill\left(\put(5,-15){\epsfbox{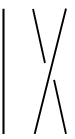}}\hspace{1cm}
 -\put(5,-15){\epsfbox{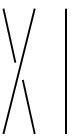}}\hspace{1cm}
 +\put(5,-15){\epsfbox{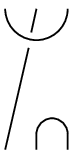}}\hspace{1cm}
 -\put(5,-15){\epsfbox{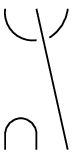}}\hspace{1cm} 
 +\put(5,-15){\epsfbox{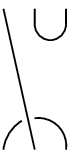}}\hspace{1cm}
 -\put(5,-15){\epsfbox{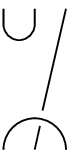}}\hspace{1cm}\right) \\  \\
+{W}^{-12}({W}^{4}-1)\\
\hfill\left(\put(5,-15){\epsfbox{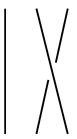}}\hspace{1cm}
 -\put(5,-15){\epsfbox{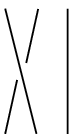}}\hspace{1cm}
 +\put(5,-15){\epsfbox{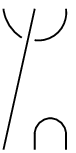}}\hspace{1cm}
 -\put(5,-15){\epsfbox{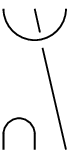}}\hspace{1cm}
 +\put(5,-15){\epsfbox{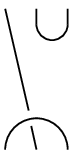}}\hspace{1cm}
 -\put(5,-15){\epsfbox{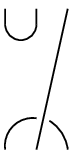}}\hspace{1cm}\right)\\ \\ 
+W^{-27}({W}^{4}-1)^{2}(1-{W}^{28}+{W}^{24}+3{W}^{16}+{W}^{12}+2{W}^{8})\\
\hfill\left(\put(5,-15){\epsfbox{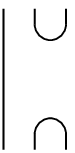}}\hspace{1cm}
 -\put(5,-15){\epsfbox{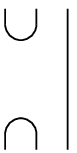}}\hspace{1cm}\right)\\ \\
\hline 
\end{array}$$
\caption{Relation ``skein'' pour $Z_{\spin_7}$}\label{skein}
\end{figure}

\begin{propri}
L'invariant d'entrelacs non orient{\'e}s $Z_{\spin_7}$ {\`a} valeurs dans $R$ poss{\`e}de
les propri{\'e}t{\'e}s suivantes~:
\begin{equation}\label{sk D}
Z_{\spin_7}\left(\put(0,-10){\epsfbox{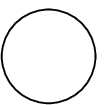}}\hspace{1cm}\amalg L\right) =
\frac{(1+W^4)(1+W^{12})(1+W^{20})}{W^{18}}Z_{\spin_7}(L)
\end{equation}
\begin{equation}\label{sk Torsion}
Z_{\spin_7}\left(\put(5,-15){\epsfbox{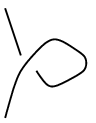}}\hspace{1cm}\right) =
W^{-21}Z_{\spin_7}\left(\put(5,-15){\epsfbox{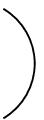}}\hspace{1cm}\right)
\end{equation}
\begin{equation}\label{sk II}
{\begin{array}{l}
Z_{\spin_7}\left(\put(5,-15){\epsfbox{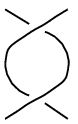}}\hspace{1cm}\right)=
W^7\,Z_{\spin_7}\left(\put(5,-15){\epsfbox{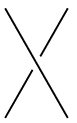}}\hspace{1cm}\right)\hfill\\\\
\hfill \hspace{1cm}
+W^4(W^{24}-1)(1+W^4)^{-1}Z_{\spin_7}
\left(\put(5,-15){\epsfbox{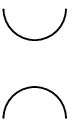}}\hspace{1cm}\right)\\\\
\hfill-W^{-15}(1+W^{-8}-W^{-12})Z_{\spin_7}
\left(\put(5,-15){\epsfbox{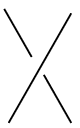}}\hspace{1cm}\right)\\\\
\hfill+W^{-2}(1+W^8-W^{12})Z_{\spin_7}
\left(\put(5,-15){\epsfbox{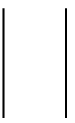}}\hspace{1cm}\right)\\
\end{array}}\end{equation}
De plus, il v{\'e}rifie la relation ``skein'' de la figure \ref{skein}.\\
On a aussi $Z_{\spin_7}(L\amalg L')=Z_{\spin_7}(L)Z_{\spin_7}(L')$ et si on note $S$
l'involution de $R$ obtenue en envoyant $\a$ sur $-\a$, alors
l'image par $Z_{\spin_7}$ de l'image miroir d'un entrelacs $L$ est donn{\'e}e par
$S_*(Z_{\spin_7}(L))$. En particulier, on obtient des relations ``skein'' encore
v{\'e}rifi{\'e}es par $Z_{\spin_7}$, {\`a} partir de celles cit{\'e}es, en changeant $W$ en
$S_*(W)=W^{-1}$ et en inversant tous les croisements.
\end{propri}

\section{Calcul de $Z_{\spin_7}$}
Remarquons que tous les coefficients intervenant dans les diff{\'e}rentes relations
``skein'' sont dans $\ZZ[W,W^{-1}]$.
\begin{theo}
Les relations ``skein'' {\'e}nonc{\'e}es permettent un calcul algorithmique de
l'invariant $Z_{\spin_7}$. En cons{\'e}quence, cet invariant d'entrelacs prend ses valeurs dans l'anneau $\ZZ[W,W^{-1}]$.
\end{theo}
\begin{dem}
Si deux segments disjoints d'un entrelacs $L$ de graphe $G(L)\in\G$ ne
s'intersectent qu'{\`a} leurs deux extr{\'e}mit{\'e}s dans le graphe $G(L)$, on appelle
{\oe}il de $G(L)$ le sous-graphe inclus dans la partie du plan, hom{\'e}omorphe {\`a} un
disque, d{\'e}limit{\'e}e par la r{\'e}union des deux segments. Un {\oe}il est minimal s'il ne
contient pas d'autres yeux que lui m{\^e}me. Si $G(L)$ ne contient pas d'yeux, alors
son invariant peut {\^e}tre calcul{\'e} par les relations (\ref{sk D}) et (\ref{sk
Torsion}).\\
Il est facile de se rendre compte que si l'int{\'e}rieur strict d'un {\oe}il minimal
n'est pas isotope {\`a} une tresse, il contient un sous-graphe pouvant {\^e}tre r{\'e}duit
par les relations (\ref{sk D}) ou (\ref{sk Torsion}). Si l'int{\'e}rieur de l'{\oe}il est
vide, par un mouvement de Reidmester de type II ou en utilisant la relation
(\ref{sk II}), on se ram{\`e}ne au calcul de l'invariant pour un graphe ayant un
nombre strictement inf{\'e}rieur de croisements. Si l'int{\'e}rieur de l'{\oe}il est une
tresse obtenue comme le produit $\sigma_1\ldots\sigma_n$ de $n$ tresses
{\'e}l{\'e}mentaires, par un mouvement de Reidmester de type III ou par la relation
figure \ref{skein}, il est possible de se ramener au calcul de l'invariant pour
des graphes ayant moins de croisements et pour un graphe ayant le m{\^e}me nombre de
croisements mais contenant un {\oe}il dont l'int{\'e}rieur est isotope {\`a} la tresse
$\sigma_2\ldots\sigma_n$~: 
$$\epsfbox{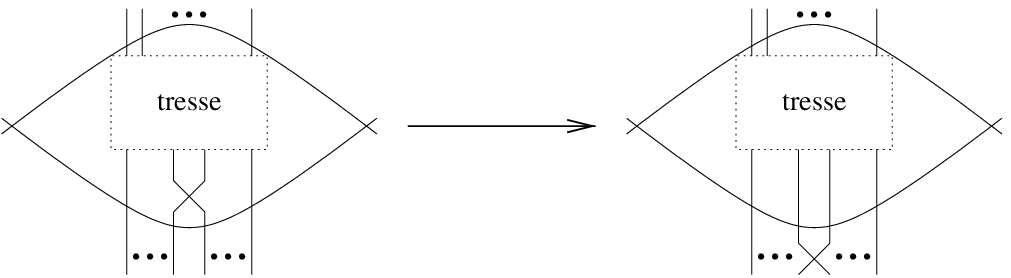}$$
En r{\'e}it{\'e}rant ce processus $n$ fois, on se ram{\`e}ne au
cas o{\`u} l'int{\'e}rieur de l'{\oe}il est form{\'e} de $k$ brins verticaux. Ensuite,
toujours par un mouvement de Reidmester de type III ou par la relation figure
\ref{skein}, il est possible se ramener au calcul de l'invariant pour
des graphes ayant moins de croisements et pour un graphe ayant le m{\^e}me nombre de
croisements mais contenant un {\oe}il dont l'int{\'e}rieur est compos{\'e} de $k-1$ brins
verticaux~: 
$$\epsfbox{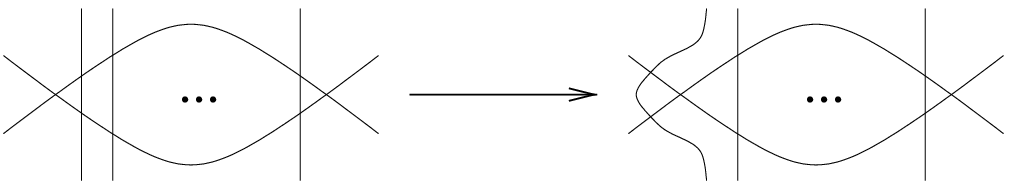}$$
En r{\'e}p{\'e}tant $k$ fois ce processus, on se ram{\`e}ne au cas o{\`u} le graphe
contient un {\oe}il vide et donc {\`a} calculer l'invariant de graphes ayant un
nombre strictement inf{\'e}rieur de croisements.
\end{dem}

\section{Exemples}\label{exemples}
En fait, les relations (\ref{sk D}) {\`a} (\ref{sk II}) suffisent pour calculer les
invariants sur des entrelacs assez simples~:

$$\Delta=Z_{\spin_7}\left(\put(0,-10){\epsfbox{B05.eps}}\hspace{1cm}\right) =
\frac{(1+W^4)(1+W^{12})(1+W^{20})}{W^{18}}$$

$$Z_{\spin_7}\left(\put(0,-8){\epsfbox{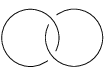}}\hspace{1.1cm}\right)=
\frac\Delta{W^{24}}(1+W^{24})(1+W^{16})(1+W^8)$$

$$Z_{\spin_7}\left(\put(0,-14){\epsfbox{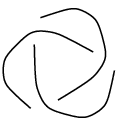}}\hspace{1.2cm}\right)=\Delta\left(
W^{27}+W^{19}-W^{15}+W^{11}-W^{7}+2W^{3}-2W^{-1}+2W^{-5}\right.$$
$$\hspace{1cm} -2W^{-9}+2W^{-13}
-2W^{-17}+2W^{-21}-2W^{-25}+W^{-29}-2W^{-33}+W^{-37}$$
$$\left.\hspace{1cm} -W^{-41}+W^{-45}\right)$$

$$Z_{\spin_7}\left(\put(0,-14){\epsfbox{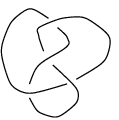}}\hspace{1.1cm}\right)=\Delta\left(
W^{48}-W^{44}+2W^{40}-3W^{36}+3W^{32}-4W^{28}+6W^{24}\right.$$
$$\hspace{1cm} -6W^{20}+7W^{16}
-8W^{12}+8W^8-9W^4+9-9W^{-4}+8W^{-8}-8W^{-12}+7W^{-16}$$
$$\left.\hspace{1cm} -6W^{-20}+6W^{-24}
-4W^{-28}+3W^{-32}-3W^{-36}+2W^{-40}-W^{-44}+W^{-48}\right)$$


\vfill

\end{document}